# ASYMPTOTIC RESULTS WITH GENERALIZED ESTIMATING EQUATIONS FOR LONGITUDINAL DATA

BY R. M. BALAN[1] AND I. SCHIOPU-KRATINA

*University of Ottawa and Statistics Canada*

We consider the marginal models of Liang and Zeger [*Biometrika* **73** (1986) 13–22] for the analysis of longitudinal data and we develop a theory of statistical inference for such models. We prove the existence, weak consistency and asymptotic normality of a sequence of estimators defined as roots of pseudo-likelihood equations.

**1. Introduction.** Longitudinal data sets arise in biostatistics and lifetime testing problems when the responses of the individuals are recorded repeatedly over a period of time. By controlling for individual differences, longitudinal studies are well-suited to measure change over time. On the other hand, they require the use of special statistical techniques because the responses on the same individual tend to be strongly correlated. In a seminal paper Liang and Zeger (1986) proposed the use of generalized linear models (GLM) for the analysis of longitudinal data.

In a cross-sectional study, a GLM is used when there are reasons to believe that each response $y_i$ depends on an observable vector $\mathbf{x}_i$ of covariates [see the monograph of McCullagh and Nelder (1989)]. Typically this dependence is specified by an unknown parameter $\beta$ and a link function $\mu$ via the relationship $\mu_i(\beta) = \mu(\mathbf{x}_i^T \beta)$, where $\mu_i(\beta)$ is the mean of $y_i$. For one-dimensional observations, the maximum quasi-likelihood estimator $\hat{\beta}_n$ is defined as the solution of the equation

$$\text{(1)} \qquad \sum_{i=1}^n \dot{\mu}_i(\beta) v_i(\beta)^{-1} (y_i - \mu_i(\beta)) = 0,$$

where $\dot{\mu}_i$ is the derivative of $\mu_i$ and $v_i(\beta)$ is the variance of $y_i$. Note that this equation simplifies considerably if we assume that $v_i(\beta) = \phi_i \dot{\mu}(\mathbf{x}_i^T \beta)$, with a

Received May 2003; revised March 2004.
[1]Supported by the Natural Sciences and Engineering Research Council of Canada.
*AMS 2000 subject classifications.* Primary 62F12; secondary 62J12.
*Key words and phrases.* Generalized estimating equations, generalized linear model, consistency, asymptotic normality.

nuisance scale parameter $\phi_i$. In fact (1) is a genuine likelihood equation if the $y_i$'s are independent with densities $c(y_i, \phi_i) \exp\{\phi_i^{-1}[(\mathbf{x}_i^T \beta)y_i - b(\mathbf{x}_i^T \beta)]\}$.

In a longitudinal study, the components of an observation $\mathbf{y}_i = (y_{i1}, \ldots, y_{im})^T$ represent repeated measurements at times $1, \ldots, m$ for subject $i$. The approach proposed by Liang and Zeger is to impose the usual assumptions of a GLM only for the marginal scalar observations $y_{ij}$ and the $p$-dimensional design vectors $\mathbf{x}_{ij}$. If the correlation matrices within individuals are known (but the entire likelihood is not specified), then the $m$-dimensional version of (1) becomes a *generalized estimating equation* (GEE).

In this article we prove the existence, weak consistency and asymptotic normality of a sequence of estimators, defined as solutions (roots) of *pseudo-likelihood equations* [see Shao (1999), page 315]. We work within a nonparametric set-up similar to that of Liang and Zeger and build upon the impressive work of Xie and Yang (2003).

Our approach differs from that of Liang and Zeger (1986), Xie and Yang (2003) and Schiopu-Kratina (2003) in the treatment of the correlation structure of the data recorded for the same individual across time. As in Rao (1998), we first obtain a sequence of preliminary consistent estimators $(\tilde{\beta}_n)_n$ of the main parameter $\beta_0$ (under the "working independence assumption"), which we use to consistently estimate the average of the true individual correlations. We then create the pseudo-likelihood equations whose solutions provide our final sequence of consistent estimators of the main parameter. In practice, the analyst would first use numerical approximation methods (like the Newton–Raphson method) to solve a simple estimating equation, where each individual correlation matrix is the identity matrix. The next step would be to solve for $\beta$ in the pseudo-likelihood equation, in which all the quantities can be calculated from the data. This approach eliminates the need to introduce nuisance parameters or to guess at the correlation structures, and thus avoids some of the problems associated with these methods [see pages 112 and 113 of Fahrmeir and Tutz (1994)]. We note that the assumptions that we require for this two-step procedure [our conditions $(\widetilde{\text{AH}})$, $(\widetilde{\text{I}}_w)$, $(\widetilde{\text{C}}_w)$] are only slightly more stringent than those of Xie and Yang (2003). They reduce to conditions related to the "working independence assumption" when the average of the true correlation matrices is asymptotically nonsingular [our hypothesis (H)].

As in Lai, Robbins and Wei (1979), where the linear model is treated, we relax the assumption of independence between subjects and consider residuals which form a martingale difference sequence. Thus our results are more general than results published so far, for example, Xie and Yang (2003) for GEE, and Shao (1992) for GLM.

Since a GEE is not a derivative, most of the technical difficulties surface when proving the existence of roots of such general estimating equations.

Two distinct methods have been developed to deal with this problem. One gives a local solution of the GEE and relies on the classical proof of the inverse function theorem [Yuan and Jennrich (1998) and Schiopu-Kratina (2003)]. The other method, which uses a result from topology, was first brought into this context by Chen, Hu and Ying (1999) and was extensively used by Xie and Yang (2003) in their proof of consistency. We adopt this second method, which facilitates a comparison of our results to those of Xie and Yang (2003) and incorporates the inference results for GLM contained in the seminal work of Fahrmeir and Kaufmann (1985).

This article is organized as follows. Section 2 is dedicated to the existence and weak consistency of a sequence of estimators of the main parameter. To accommodate the estimation of the average of the correlation matrices in the martingale set-up, we require two conditions: (C1) is a boundedness condition on the $(2+\delta)$-moments of the normalized residuals, whereas (C2) is a consistency condition on the normalized conditional covariance matrix. In this context we use the martingale strong law of large numbers of Kaufmann (1987). Section 3 presents the asymptotic normality of our estimators. This is obtained under slightly stronger conditions than those of Xie and Yang (2003), by applying the classical martingale central limit theorem [see Hall and Heyde (1980)]. For ease of exposition, we have placed the more technical proofs in the Appendix.

We introduce first some matrix notation [see Schott (1997)]. If $\mathbf{A}$ is a $p \times p$ matrix, we will denote with $\|\mathbf{A}\|$ its spectral norm, with $\det(\mathbf{A})$ its determinant and with $\text{tr}(\mathbf{A})$ its trace. If $\mathbf{A}$ is a symmetric matrix, we denote by $\lambda_{\min}(\mathbf{A})[\lambda_{\max}(\mathbf{A})]$ its minimum (maximum) eigenvalue. For any matrix $\mathbf{A}$, $\|\mathbf{A}\| = \{\lambda_{\max}(\mathbf{A}^T \mathbf{A})\}^{1/2}$. For a $p$-dimensional vector $\mathbf{x}$, we use the Euclidean norm $\|\mathbf{x}\| = (\mathbf{x}^T \mathbf{x})^{1/2} = \text{tr}(\mathbf{x}\mathbf{x}^T)^{1/2}$. We let $\mathbf{A}^{1/2}$ be the symmetric square root of a positive definite matrix $\mathbf{A}$ and $\mathbf{A}^{-1/2} = (\mathbf{A}^{1/2})^{-1}$. Finally, we use the matrix notation $\mathbf{A} \leq \mathbf{B}$ if $\lambda^T \mathbf{A} \lambda \leq \lambda^T \mathbf{B} \lambda$ for any $p$-dimensional vector $\lambda$.

Throughout this article we will assume that the number of longitudinal observations on each individual is fixed and equal to $m$. More precisely, we will denote with $\mathbf{y}_i := (y_{i1}, \ldots, y_{im})'$, $i \leq n$, a longitudinal data set consisting of $n$ respondents, where the components of $\mathbf{y}_i$ represent measurements at different times on subject $i$. The observations $y_{ij}$ are recorded along with a corresponding $p$-dimensional vector $\mathbf{x}_{ij}$ of covariates and the marginal expectations and variances are specified in terms of the regression parameter $\beta$ through $\theta_{ij} = \mathbf{x}_{ij}^T \beta$ as follows:

(2) $\quad \mu_{ij}(\beta) := E_\beta(y_{ij}) = \mu(\theta_{ij}), \qquad \sigma_{ij}^2(\beta) := \text{Var}_\beta(y_{ij}) = \dot{\mu}(\theta_{ij}),$

where $\mu$ is a continuously differentiable link function with $\dot{\mu} > 0$, that is, we consider only canonical link functions.

Here are the most commonly used such link functions:

1. In the linear regression, $\mu(y) = y$.
2. In the log regression for count data, $\mu(y) = \exp(y)$.
3. In the logistic regression for binary data, $\mu(y) = \exp(y)/[1+\exp(y)]$.
4. In the probit regression for binary data, $\mu(y) = \Phi(y)$, where $\Phi$ is the standard normal distribution function; we have $\dot{\Phi}(y) = (2\pi)^{-1/2} \exp(-y^2/2)$.

In the sequel the unknown parameter $\beta$ lies in an open set $\mathcal{B} \subseteq R^p$ and $\beta_0$ is the true value of this parameter. We normally drop the parameter $\beta_0$ to avoid cumbersome notation.

Let $\mu_i(\beta) = (\mu_{i1}(\beta), \ldots, \mu_{im}(\beta))^T$, $\mathbf{A}_i(\beta) = \mathrm{diag}(\sigma_{i1}^2(\beta), \ldots, \sigma_{im}^2(\beta))$ and $\Sigma_i(\beta) := \mathrm{Cov}_\beta(\mathbf{y}_i)$. Note that $\Sigma_i = \mathbf{A}_i^{1/2} \bar{\mathbf{R}}_i \mathbf{A}_i^{1/2}$, where $\bar{\mathbf{R}}_i$ is the true correlation matrix of $\mathbf{y}_i$ at $\beta_0$. Let $\mathbf{X}_i = (\mathbf{x}_{i1}, \ldots, \mathbf{x}_{im})^T$.

We consider the sequence $\varepsilon_i(\beta) = (\varepsilon_{i1}(\beta), \ldots, \varepsilon_{im}(\beta))^T$ with $\varepsilon_{ij}(\beta) = \mathbf{y}_{ij} - \mu_{ij}(\beta)$, and we assume that the residuals $(\varepsilon_i)_{i \geq 1}$ form a martingale difference sequence, that is,

$$E(\varepsilon_i | \mathcal{F}_{i-1}) = 0 \qquad \text{for all } i \geq 1,$$

where $\mathcal{F}_i$ is the minimal $\sigma$-field with respect to which $\varepsilon_1, \ldots, \varepsilon_i$ are measurable. This is a natural generalization of the case of independent observations.

Finally, to avoid keeping track of various constants, we agree to denote with $C$ a generic constant which does not depend on $n$, but is different from case to case.

**2. Asymptotic existence and consistency.** We consider the generalized estimating equations (GEE) of Xie and Yang (2003) in the case when the "working" correlation matrices are $\mathbf{R}_i^{\mathrm{indep}} = \mathbf{I}$ for all $i$. This is also known as the "working independence" case, the word "independence" referring to the observations on the same individual. Let $(\tilde{\beta}_n)_n$ be a sequence of estimators such that

(3) $$P(\mathbf{g}_n^{\mathrm{indep}}(\tilde{\beta}_n) = 0) \to 1 \quad \text{and} \quad \tilde{\beta}_n \xrightarrow{P} \beta_0,$$

where $\mathbf{g}_n^{\mathrm{indep}}(\beta) = \sum_{i=1}^n \mathbf{X}_i^T \varepsilon_i(\beta)$ is the "working independence" GEE.

The following quantities have been used extensively in the work of Xie and Yang (2003) and play an important role in the conditions for the existence and consistency of $\tilde{\beta}_n$:

$$\mathbf{H}_n^{\mathrm{indep}} = \sum_{i=1}^n \mathbf{X}_i^T \mathbf{A}_i \mathbf{X}_i, \qquad \pi_n^{\mathrm{indep}} := \frac{\max_{i \leq n} \lambda_{\max}((\mathbf{R}_i^{\mathrm{indep}})^{-1})}{\min_{i \leq n} \lambda_{\min}((\mathbf{R}_i^{\mathrm{indep}})^{-1})} = 1,$$

$$\tilde{\tau}_n^{\mathrm{indep}} := m \max_{i \leq n} \lambda_{\max}((\mathbf{R}_i^{\mathrm{indep}})^{-1}) = m,$$

$$(\gamma_n^{(0)})^{\mathrm{indep}} := \max_{i \leq n, j \leq m} \mathbf{x}_{ij}^T (\mathbf{H}_n^{\mathrm{indep}})^{-1} \mathbf{x}_{ij}.$$

We will also use the following maxima:

$$k_n^{[2]}(\beta) = \max_{i \leq n} \max_{j \leq m} \left| \frac{\ddot{\mu}(\theta_{ij})}{\dot{\mu}(\theta_{ij})} \right|, \qquad k_n^{[3]}(\beta) = \max_{i \leq n} \max_{j \leq m} \left| \frac{\mu^{(3)}(\theta_{ij})}{\dot{\mu}(\theta_{ij})} \right|.$$

The fact that the residuals $(\varepsilon_i)_{i \geq 1}$ form a martingale difference sequence does not change the proofs of Theorem 2 and Theorem A.1(ii) of Xie and Yang (2003). Following their work, we conclude that the sufficient conditions for the existence of a sequence $(\tilde{\beta}_n)_n$ with the desired property (3) are:

(AH)$^{\text{indep}}$ for any $r > 0$, $k_n^{[l],\text{indep}} = \sup_{\beta \in B_n^{\text{indep}}(r)} k_n^{[l]}(\beta)$, $l = 2, 3$, are bounded,

$(\text{I}_w^*)^{\text{indep}}$ $\lambda_{\min}(\mathbf{H}_n^{\text{indep}}) \to \infty$,

$(\text{C}_w^*)^{\text{indep}}$ $n^{1/2}(\gamma_n^{(0)})^{\text{indep}} \to 0$,

where $B_n^{\text{indep}}(r) := \{\beta; \|(\mathbf{H}_n^{\text{indep}})^{1/2}(\beta - \beta_0)\| \leq m^{1/2}r\}$. We denote by (C)$^{\text{indep}}$ the set of conditions (AH)$^{\text{indep}}$, $(\text{I}_w^*)^{\text{indep}}$, $(\text{C}_w^*)^{\text{indep}}$.

It turns out that, in practice, the analyst will have to verify conditions similar to (C)$^{\text{indep}}$ in order to produce the estimators that we propose (see Remark 5). All the classical examples corresponding to our link functions 1–4 are within the scope of our theory. We present below two new examples.

EXAMPLE 1. Suppose that $p = 2$. Let $\mathbf{x}_{ij} = (a_{ij}, b_{ij})^T$, $u_n = \sum_{i \leq n, j \leq m} \sigma_{ij}^2 a_{ij}^2$, $v_n = \sum_{i \leq n, j \leq m} \sigma_{ij}^2 b_{ij}^2$ and $w_n = \sum_{i \leq n, j \leq m} \sigma_{ij}^2 a_{ij} b_{ij}$. In this case

$$\mathbf{H}_n^{\text{indep}} = \begin{bmatrix} u_n & w_n \\ w_n & v_n \end{bmatrix},$$

$\lambda_{\max}(\mathbf{H}_n^{\text{indep}}) = (u_n + v_n + d_n)/2$ and $\lambda_{\min}(\mathbf{H}_n^{\text{indep}}) = (u_n + v_n - d_n)/2$, with $d_n := \sqrt{(u_n - v_n)^2 + 4w_n^2}$. Note that $w_n = \sqrt{u_n v_n} \cos \theta_n$ for $\theta_n \in [0, \pi]$ and $\det(\mathbf{H}_n^{\text{indep}}) = u_n v_n \sin^2 \theta_n$ [see also page 79 of McCullagh and Nelder (1989)]. Suppose that

$$\alpha := \liminf_{n \to \infty} \sin^2 \theta_n > 0.$$

Since

$$\frac{1}{\lambda_{\min}(\mathbf{H}_n^{\text{indep}})} = \frac{\lambda_{\max}(\mathbf{H}_n^{\text{indep}})}{\det(\mathbf{H}_n^{\text{indep}})} = \frac{u_n + v_n + d_n}{u_n v_n \sin^2 \theta_n},$$

one can show that condition $(\text{I}_w^*)^{\text{indep}}$ is equivalent to $\min(u_n, v_n) \to \infty$. On the other hand,

$$\mathbf{x}_{ij}^T (\mathbf{H}_n^{\text{indep}})^{-1} \mathbf{x}_{ij} = \frac{a_{ij}^2}{u_n} - 2w_n \frac{a_{ij} b_{ij}}{u_n v_n} + \frac{b_{ij}^2}{v_n} \leq \left( \frac{a_{ij}}{u_n^{1/2}} + \frac{b_{ij}}{v_n^{1/2}} \right)^2.$$

Condition $(\text{C}_w^*)^{\text{indep}}$ holds if $n^{1/2} \max_{i \leq n, j \leq m} (u_n^{-1/2} a_{ij} + v_n^{-1/2} b_{ij})^2 \to 0$.

EXAMPLE 2. The case of a single covariate with $p$ different levels [one-way ANOVA; see also Example 3.13 of Shao (1999)] is usually treated by identifying each of these levels with one of the $p$-dimensional vectors $\mathbf{e}_1, \ldots, \mathbf{e}_p$, where $\mathbf{e}_k$ has the $k$th component 1 and all the other components 0. We can say that $\mathbf{x}_{ij} \in \{\mathbf{e}_1, \ldots, \mathbf{e}_p\}$ for all $i \leq n, j \leq m$. In this case, $\mathbf{H}_n^{\text{indep}}$ is a diagonal matrix. More precisely,

$$\mathbf{H}_n^{\text{indep}} = \sum_{k=1}^{p} \nu_n^{(k)} \mathbf{e}_k \mathbf{e}_k^T,$$

where $\nu_n^{(k)} = \sum_{i \leq n, j \leq m; \mathbf{x}_{ij} = \mathbf{e}_k} \sigma_{ij}^2$. Let $\nu_n = \min_{k \leq p} \nu_n^{(k)}$. Condition $(\text{I}_w^*)^{\text{indep}}$ is equivalent to $\nu_n \to \infty$ and condition $(\text{C}_w^*)^{\text{indep}}$ is equivalent to $n^{1/2} \nu_n^{-1} \to 0$.

The method introduced by Liang and Zeger (1986) and developed recently in Xie and Yang (2003) relies heavily on the "working" correlation matrices $\mathbf{R}_i(\alpha)$ which are chosen arbitrarily by the statistician (possibly containing a nuisance parameter $\alpha$) and are expected to be good approximations of the unknown true correlation matrices $\bar{\mathbf{R}}_i$.

In the present paper, we consider an alternative approach in which at each step $n$, the "working" correlation matrices $\mathbf{R}_i(\alpha), i \leq n$, are replaced by the random matrix

$$\widetilde{\mathcal{R}}_n := \frac{1}{n} \sum_{i=1}^{n} \mathbf{A}_i(\tilde{\beta}_n)^{-1/2} \varepsilon_i(\tilde{\beta}_n) \varepsilon_i(\tilde{\beta}_n)^T \mathbf{A}_i(\tilde{\beta}_n)^{-1/2}$$

which depends only on the data set and is shown to be a (possibly biased) consistent estimator of the average of the true correlation matrices

$$\bar{\bar{\mathbf{R}}}_n := \frac{1}{n} \sum_{i=1}^{n} \bar{\mathbf{R}}_i.$$

The consistency of $\widetilde{\mathcal{R}}_n$ is obtained under the following two conditions imposed on the (normalized) residuals $\mathbf{y}_i^* = \mathbf{A}_i^{-1/2} \varepsilon_i$, with $E(\mathbf{y}_i^* \mathbf{y}_i^{*T}) = \bar{\mathbf{R}}_i$:

(C1) there exists a $\delta \in (0, 2]$ such that $\sup_{i \geq 1} E(\|\mathbf{y}_i^*\|^{2+\delta}) < \infty$,
(C2) $\frac{1}{n} \sum_{i=1}^{n} \mathcal{V}_i \xrightarrow{P} 0$, where $\mathcal{V}_i = E(\mathbf{y}_i^* \mathbf{y}_i^{*T} | \mathcal{F}_{i-1}) - \bar{\mathbf{R}}_i$.

REMARK 1. Condition (C1) is a bounded moment requirement which is usually needed for verifying the conditions of a martingale limit theorem, while condition (C2) is satisfied if the observations are independent. Condition (C2) is in fact a requirement on the (normalized) conditional covariance matrix $\mathbf{V}_n = \sum_{i=1}^{n} E(\mathbf{y}_i^* \mathbf{y}_i^{*T} | \mathcal{F}_{i-1})$. More precisely, if the following hypothesis holds true:

(H) there exists a constant $C > 0$ such that $\lambda_{\min}(\bar{\bar{\mathbf{R}}}_n) \geq C$ for all $n$,

then condition (C2) is equivalent to $\bar{\bar{\mathbf{R}}}_n^{-1/2}(\mathbf{V}_n/n)\bar{\bar{\mathbf{R}}}_n^{-1/2} - \mathbf{I} \xrightarrow{P} 0$ [which is similar to (3.1) of Hall and Heyde (1980) or (4.2) of Shao (1992)]. Note that (H) is implied by the following stronger hypothesis, which is needed in Section 3:

(H′) There exists a constant $C > 0$ such that $\lambda_{\min}(\bar{\mathbf{R}}_i) \geq C$ for all $i$.

Hypothesis (H′) is satisfied if $\bar{\mathbf{R}}_i = \bar{\mathbf{R}}$ for all $i$, where $\bar{\mathbf{R}}$ is nonsingular.

The following result is essential for all our developments.

THEOREM 1. *Let $\mathbf{R}_n = E(\widetilde{\mathcal{R}}_n)$. Under conditions (C)$^{\text{indep}}$, (C1) and (C2), we have*

$$\widetilde{\mathcal{R}}_n - \mathbf{R}_n \xrightarrow{L^1} 0 \qquad (elementwise).$$

*If the convergence in condition (C2) is almost sure, then $\widetilde{\mathcal{R}}_n - \mathbf{R}_n \xrightarrow{a.s.} 0$ (elementwise). The same conclusion holds if $\mathbf{R}_n$ is replaced by $\bar{\bar{\mathbf{R}}}_n$.*

PROOF. Let $\widehat{\mathcal{R}}_n = n^{-1}\sum_{i=1}^n \mathbf{A}_i^{-1/2}\varepsilon_i\varepsilon_i^T \mathbf{A}_i^{-1/2}$ and note that $E(\widehat{\mathcal{R}}_n) = \bar{\bar{\mathbf{R}}}_n$. Our result will be a consequence of the following two propositions, whose proofs are given in Appendix A. □

PROPOSITION 1. *Under conditions (C1) and (C2), we have*

$$\widehat{\mathcal{R}}_n - \bar{\bar{\mathbf{R}}}_n \xrightarrow{L^1} 0 \qquad (elementwise).$$

PROPOSITION 2. *Under conditions (C)$^{\text{indep}}$, (C1) and (C2), we have*

$$\widetilde{\mathcal{R}}_n - \widehat{\mathcal{R}}_n \xrightarrow{L^1} 0 \qquad (elementwise).$$

In what follows we will assume that the inverse of the (nonnegative definite) random matrix $\widetilde{\mathcal{R}}_n$ exists with probability 1, for every $n$. We consider the following *pseudo-likelihood* equation:

$$(4) \qquad \sum_{i=1}^n \mathbf{D}_i(\beta)^T \widetilde{\mathbf{V}}_{i,n}(\beta)^{-1}\varepsilon_i(\beta) = 0,$$

where $\mathbf{D}_i(\beta) = \mathbf{A}_i(\beta)\mathbf{X}_i$ and $\widetilde{\mathbf{V}}_{i,n}(\beta) := \mathbf{A}_i(\beta)^{1/2}\widetilde{\mathcal{R}}_n \mathbf{A}_i(\beta)^{1/2}$. Note that (4) can be written as

$$\tilde{\mathbf{g}}_n(\beta) := \sum_{i=1}^n \mathbf{X}_i^T \mathbf{A}_i(\beta)^{1/2}\widetilde{\mathcal{R}}_n^{-1}\mathbf{A}_i(\beta)^{-1/2}\varepsilon_i(\beta) = 0.$$

We consider also the estimating function

$$\mathbf{g}_n(\beta) = \sum_{i=1}^{n} \mathbf{X}_i^T \mathbf{A}_i(\beta)^{1/2} \mathbf{R}_n^{-1} \mathbf{A}_i(\beta)^{-1/2} \varepsilon_i(\beta).$$

Note that $\mathbf{M}_n := \text{Cov}(\mathbf{g}_n) = \sum_{i=1}^{n} \mathbf{X}_i^T \mathbf{A}_i^{1/2} \mathbf{R}_n^{-1} \bar{\mathbf{R}}_i \mathbf{R}_n^{-1} \mathbf{A}_i^{1/2} \mathbf{X}_i$.

As in Xie and Yang (2003), we introduce the following quantities:

$$\mathbf{H}_n := \sum_{i=1}^{n} \mathbf{X}_i^T \mathbf{A}_i^{1/2} \mathbf{R}_n^{-1} \mathbf{A}_i^{1/2} \mathbf{X}_i, \qquad \pi_n := \frac{\lambda_{\max}(\mathbf{R}_n^{-1})}{\lambda_{\min}(\mathbf{R}_n^{-1})},$$

$$\tilde{\tau}_n := m \lambda_{\max}(\mathbf{R}_n^{-1}),$$

$$\gamma_n^{(0)} := \max_{i=1,\ldots,n} \max_{j=1,\ldots,m} (\mathbf{x}_{ij}^T \mathbf{H}_n^{-1} \mathbf{x}_{ij}), \qquad \tilde{\gamma}_n = \tilde{\tau}_n \gamma_n^{(0)}.$$

REMARK 2. A few comments about $\tilde{\tau}_n$ are worth mentioning. First, $\mathbf{M}_n \leq \tau_n \mathbf{H}_n$, where $\tau_n := \max_{i \leq n} \lambda_{\max}(\mathbf{R}_n^{-1} \bar{\mathbf{R}}_i) \leq \tilde{\tau}_n$. Also, since $r_{jk}^{(n)} - \bar{\bar{r}}_{jk}^{(n)} \to 0$ and $|\bar{\bar{r}}_{jk}^{(n)}| \leq 1$, we can assume that $|r_{jk}^{(n)}| \leq 2$, for $n$ large enough (here $r_{jk}^{(n)}, \bar{\bar{r}}_{jk}^{(n)}$ are the elements of the matrices $\mathbf{R}_n$, resp. $\bar{\bar{\mathbf{R}}}_n$). Therefore $\tilde{\tau}_n \geq 1/2$. The reason why we prefer to work with $\tilde{\tau}_n$ instead of $\tau_n$ will become apparent in the proof of Proposition 3 (given in Appendix A.2). Another reason is, of course, the fact that $\tilde{\tau}_n$ does not depend on the unknown matrices $\bar{\mathbf{R}}_i$.

Our approach requires a slight modification of the conditions introduced by Xie and Yang (2003) to accommodate the use of $\tilde{\tau}_n$ instead of $\tau_n$. Let $\widetilde{B}_n(r) := \{\beta; \|\mathbf{H}_n^{1/2}(\beta - \beta_0)\| \leq (\tilde{\tau}_n)^{1/2} r\}$. Our conditions are:

$(\widetilde{\text{AH}})$ for any $r > 0$, $\tilde{k}_n^{[l]} = \sup_{\beta \in \widetilde{B}_n(r)} k_n^{[l]}(\beta)$, $l = 2, 3$, are bounded,

$(\widetilde{\text{I}}_w)$ $(\tilde{\tau}_n)^{-1} \lambda_{\min}(\mathbf{H}_n) \to \infty$,

$(\widetilde{\text{C}}_w)$ $(\pi_n)^2 \tilde{\gamma}_n \to 0$, and $n^{1/2} \pi_n \tilde{\gamma}_n \to 0$.

REMARK 3. Note that $(\widetilde{\text{I}}_w)$ implies $(\text{I}_w^*)^{\text{indep}}$, which implies $\lambda_{\min}(\mathbf{H}_n) \to \infty$. This follows from the inequalities

$$\frac{1}{2m} \mathbf{H}_n^{\text{indep}} \leq \lambda_{\min}(\mathbf{R}_n^{-1}) \cdot \mathbf{H}_n^{\text{indep}} \leq \mathbf{H}_n \leq \lambda_{\max}(\mathbf{R}_n^{-1}) \cdot \mathbf{H}_n^{\text{indep}} = \frac{\tilde{\tau}_n}{m} \mathbf{H}_n^{\text{indep}}.$$

REMARK 4. Our conditions depend on the matrix $\mathbf{R}_n$, which cannot be written in a closed form. Since $\widetilde{\mathcal{R}}_n - \mathbf{R}_n \xrightarrow{P} 0$, it is desirable to express our conditions in terms of the matrix $\widetilde{\mathcal{R}}_n$. In practice, if the sample size is large enough, one may choose to verify conditions $(\widetilde{\text{AH}})$, $(\widetilde{\text{I}}_w)$, $(\widetilde{\text{C}}_w)$ by using $\widetilde{\mathcal{R}}_n$ (instead of $\mathbf{R}_n$) in the definitions of $\mathbf{H}_n, \pi_n, \tilde{\gamma}_n$.

REMARK 5. If we suppose that hypothesis (H) holds, then for $n$ large

$$\frac{C}{2} \leq \lambda_{\min}(\mathbf{R}_n) \leq \lambda_{\max}(\mathbf{R}_n) \leq 2m.$$

In this case $(\tilde{\tau}_n)_n$ and $(\pi_n)_n$ are bounded, $C'(\gamma_n^{(0)})^{\text{indep}} \leq \gamma_n^{(0)} \leq C(\gamma_n^{(0)})^{\text{indep}}$, and for every $r > 0$ there exists $r' > 0$ such that $\widetilde{B}_n(r) \subseteq B_n^{\text{indep}}(r')$. Therefore, conditions $(\widetilde{\text{AH}})$, $(\widetilde{\text{I}}_w)$, $(\widetilde{\text{C}}_w)$ are equivalent to $(\text{AH})^{\text{indep}}$, $(\text{I}_w^*)^{\text{indep}}$, $(\text{C}_w^*)^{\text{indep}}$, respectively. In order to verify (H), it is sufficient to check that there exists a constant $C > 0$ such that

$$\det(\widetilde{\mathcal{R}}_n) \geq C \qquad \text{for all } n \text{ a.s.}$$

under the hypothesis of Theorem 1.

We need to consider the derivatives

$$\widetilde{\mathcal{D}}_n(\beta) := -\frac{\partial \tilde{\mathbf{g}}_n(\beta)}{\partial \beta^T}, \qquad \mathcal{D}_n(\beta) := -\frac{\partial \mathbf{g}_n(\beta)}{\partial \beta^T}.$$

The next theorem is a modified version of Theorem A.2, respectively, Theorem A.1(ii) of Xie and Yang (2003).

THEOREM 2. *Under conditions $(\widetilde{\text{AH}})$ and $(\widetilde{\text{C}}_w)$:*

(i) *for every $r > 0$*

$$\sup_{\beta \in \widetilde{B}_n(r)} \|\mathbf{H}_n^{-1/2} \mathcal{D}_n(\beta) \mathbf{H}_n^{-1/2} - \mathbf{I}\| \xrightarrow{P} 0;$$

(ii) *there exists $c_0 > 0$ such that for every $r > 0$*

$$P(\mathcal{D}_n(\beta) \geq c_0 \mathbf{H}_n \text{ for all } \beta \in \widetilde{B}_n(r)) \to 1.$$

PROOF. (i) The first two terms produced by the decomposition $\mathcal{D}_n(\beta) = \mathbf{H}_n(\beta) + \mathbf{B}_n(\beta) + \mathcal{E}_n(\beta)$ are shown to be bounded by $\pi_n^2 \tilde{\gamma}_n$, whereas the third term is bounded in $L^2$ by $\sqrt{n} \pi_n \tilde{\gamma}_n$. [Here $\mathbf{H}_n(\beta), \mathbf{B}_n(\beta), \mathcal{E}_n(\beta)$ have the same expressions as those given in Xie and Yang (2003) with $\mathbf{R}_i(\alpha)$, $i \leq n$, replaced by $\mathbf{R}_n$.] The arguments are essentially the same as those used in Lemmas A.1(ii), A.2(ii) and A.3(ii) of Xie and Yang (2003). The fact that we are replacing the "working" correlation matrices $\mathbf{R}_i(\alpha)$, $i = 1, \ldots, n$, with the matrix $\mathbf{R}_n$ and we assume that $(\varepsilon_i)_{i \geq 1}$ is a martingale difference sequence does not influence the proof. Finally we note that (ii) is a consequence of (i). □

The next two results are intermediate steps that are used in the proof of our main result. Their proofs are given in Appendix A.2.

PROPOSITION 3. *Suppose that the conditions of Theorem 1 hold. Then*
$$(\tilde{\tau}_n)^{-1/2}\mathbf{H}_n^{-1/2}(\tilde{\mathbf{g}}_n - \mathbf{g}_n) \xrightarrow{P} 0.$$

PROPOSITION 4. *Suppose that the conditions of Theorem 1 hold. Under conditions* $(\widetilde{\mathrm{AH}})$ *and* $(\widetilde{\mathrm{C}}_w)$,
$$\sup_{\beta \in \widetilde{B}_n(r)} \|\mathbf{H}_n^{-1/2}[\widetilde{\mathcal{D}}_n(\beta) - \mathcal{D}_n(\beta)]\mathbf{H}_n^{-1/2}\| \xrightarrow{P} 0.$$

The next theorem is our main result. It shows that under our slightly modified conditions $(\widetilde{\mathrm{AH}}), (\widetilde{\mathrm{I}}_w), (\widetilde{\mathrm{C}}_w)$ and the additional conditions of Theorem 1, one can obtain a solution $\hat{\beta}_n$ of the pseudo-likelihood equation $\tilde{\mathbf{g}}_n(\beta) = 0$, which is also a consistent estimator of $\beta_0$.

THEOREM 3. *Suppose that the conditions of Theorem 1 hold. Under conditions* $(\widetilde{\mathrm{AH}})$, $(\widetilde{\mathrm{I}}_w)$ *and* $(\widetilde{\mathrm{C}}_w)$, *there exists a sequence* $(\hat{\beta}_n)_n$ *of random variables such that*
$$P(\tilde{\mathbf{g}}_n(\hat{\beta}_n) = 0) \to 1 \quad and \quad \hat{\beta}_n \xrightarrow{P} \beta_0.$$

PROOF. Let $\varepsilon > 0$ be arbitrary and $r = r(\varepsilon) = \sqrt{(24p)/(c_1^2 \varepsilon)}$, where $c_1$ is a constant to be specified later. We consider the events
$$\widetilde{E}_n := \left\{ \|\mathbf{H}_n^{-1/2}\tilde{\mathbf{g}}_n\| \leq \inf_{\beta \in \partial \widetilde{B}_n(r)} \|\mathbf{H}_n^{-1/2}(\tilde{\mathbf{g}}_n(\beta) - \tilde{\mathbf{g}}_n)\| \right\},$$
$$\widetilde{\Omega}_n := \{\widetilde{\mathcal{D}}_n(\bar{\beta}) \text{ nonsingular, for all } \bar{\beta} \in \widetilde{B}_n(r)\}.$$
By Lemma A of Chen, Hu and Ying (1999), it follows that on the event $\widetilde{E}_n \cap \widetilde{\Omega}_n$, there exists $\hat{\beta}_n \in \widetilde{B}_n(r)$ such that $\tilde{g}_n(\hat{\beta}_n) = 0$. Therefore, it remains to prove that $P(\widetilde{E}_n \cap \widetilde{\Omega}_n) > 1 - \varepsilon$ for $n$ large.

By Taylor's formula and Lemma 1 of Xie and Yang (2003) we obtain that for any $\beta \in \partial \widetilde{B}_n(r)$ there exist $\bar{\beta} \in \widetilde{B}_n(r)$ and a $p \times 1$ vector $\lambda$, $\|\lambda\| = 1$ such that
$$\|\mathbf{H}_n^{-1/2}(\tilde{\mathbf{g}}_n(\beta) - \tilde{\mathbf{g}}_n)\|$$
$$\geq |\lambda^T \mathbf{H}_n^{-1/2} \widetilde{\mathcal{D}}_n(\bar{\beta}) \mathbf{H}_n^{-1/2} \lambda| \cdot r(\tilde{\tau}_n)^{1/2}$$
$$\geq \{|\lambda^T \mathbf{H}_n^{-1/2} \mathcal{D}_n(\bar{\beta}) \mathbf{H}_n^{-1/2} \lambda|$$
$$\quad - |\lambda^T \mathbf{H}_n^{-1/2} [\widetilde{\mathcal{D}}_n(\bar{\beta}) - \mathcal{D}_n(\bar{\beta})] \mathbf{H}_n^{-1/2} \lambda|\} \cdot r(\tilde{\tau}_n)^{1/2}.$$

By Theorem 2(ii) there exists $c_0 > 0$ such that
$$(5) \quad P(\lambda^T \mathbf{H}_n^{-1/2} \mathcal{D}_n(\beta) \mathbf{H}_n^{-1/2} \lambda \geq c_0$$
$$\text{for all } \beta \in \widetilde{B}_n(r), \text{ for all } \lambda, \|\lambda\| = 1) > 1 - \varepsilon/6$$

when $n$ is large. Let $c'_0 \in (0, c_0)$ be arbitrary. By Proposition 4,

$$P(|\lambda^T \mathbf{H}_n^{-1/2}[\widetilde{\mathcal{D}}_n(\beta) - \mathcal{D}_n(\beta)]\mathbf{H}_n^{-1/2}\lambda| \leq c'_0 \tag{6}$$
$$\text{for all } \beta \in \widetilde{B}_n(r), \text{ for all } \lambda) > 1 - \varepsilon/6$$

when $n$ is large. Therefore, if we put $c_1 := c_0 - c'_0$, we have

$$P\left(\inf_{\beta \in \partial \widetilde{B}_n(r)} \|\mathbf{H}_n^{-1/2}(\tilde{\mathbf{g}}_n(\beta) - \tilde{\mathbf{g}}_n)\| \geq c_1 r(\tilde{\tau}_n)^{1/2}\right) > 1 - \varepsilon/3. \tag{7}$$

From (5) and (6) we can also conclude that $P(\widetilde{\Omega}_n) > 1 - \varepsilon/3$ for $n$ large.

On the other hand, by Chebyshev's inequality and our choice of $r$, we have $P(\|\mathbf{H}_n^{-1/2}\mathbf{g}_n\| \leq c_1 r(\tilde{\tau}_n)^{1/2}/2) > 1 - \varepsilon/6$ for all $n$. By Proposition 3, $P(\|\mathbf{H}_n^{-1/2}(\tilde{\mathbf{g}}_n - \mathbf{g}_n)\| \leq c_1 r(\tilde{\tau}_n)^{1/2}/2) > 1 - \varepsilon/6$ for $n$ large. Hence

$$P(\|\mathbf{H}_n^{-1/2}\tilde{\mathbf{g}}_n\| \leq c_1 r(\tilde{\tau}_n)^{1/2}) > 1 - \varepsilon/3. \tag{8}$$

From (7) and (8) we obtain that $P(\widetilde{E}_n) > 1 - (2\varepsilon)/3$ for $n$ large. This concludes the proof of the asymptotic existence.

We proceed now with the proof of the weak consistency. Let $\delta > 0$ be arbitrary. By $(\widetilde{\mathrm{I}}_w)$ we have $\tilde{\tau}_n/\lambda_{\min}(\mathbf{H}_n) < (\delta/r)^2$ for $n$ large. We know that on the event $\widetilde{E}_n \cap \widetilde{\Omega}_n$, there exists $\hat{\beta}_n \in \widetilde{B}_n(r)$ such that $\tilde{\mathbf{g}}_n(\hat{\beta}_n) = 0$. Therefore, on this event

$$\|\hat{\beta}_n - \beta_0\| \leq \|\mathbf{H}_n^{-1/2}\| \cdot \|\mathbf{H}_n^{1/2}(\hat{\beta}_n - \beta_0)\| \leq [\lambda_{\min}(\mathbf{H}_n)]^{-1/2} \cdot (\tilde{\tau}_n)^{1/2} r < \delta$$

for $n$ large. This proves that $P(\|\hat{\beta}_n - \beta_0\| \leq \delta) > 1 - \varepsilon$ for $n$ large. □

**3. Asymptotic normality.** Let $c_n = \lambda_{\max}(\mathbf{M}_n^{-1}\mathbf{H}_n)$. In this section we will suppose that $(c_n \tilde{\tau}_n)_n$ is bounded.

THEOREM 4. *Under the conditions of Theorem 3,*

$$\mathbf{M}_n^{-1/2}\tilde{\mathbf{g}}_n = \mathbf{M}_n^{-1/2}\mathbf{H}_n(\hat{\beta}_n - \beta_0) + o_P(1).$$

PROOF. On the set $\{\tilde{\mathbf{g}}_n(\hat{\beta}_n) = 0, \hat{\beta}_n \in \widetilde{B}_n(r)\}$, we have $\tilde{\mathbf{g}}_n = \widetilde{\mathcal{D}}_n(\bar{\beta}_n)(\hat{\beta}_n - \beta_0)$ for some $\bar{\beta}_n \in \widetilde{B}_n(r)$ by Taylor's formula. Multiplication with $\mathbf{M}_n^{-1/2}$ yields

$$\mathbf{M}_n^{-1/2}\tilde{\mathbf{g}}_n = \mathbf{M}_n^{-1/2}\mathbf{H}_n^{1/2}\mathbf{A}_n\mathbf{H}_n^{1/2}(\hat{\beta}_n - \beta_0) + \mathbf{M}_n^{-1/2}\mathbf{H}_n(\hat{\beta}_n - \beta_0),$$

where $\mathbf{A}_n := \mathbf{H}_n^{-1/2}\widetilde{\mathcal{D}}_n(\bar{\beta}_n)\mathbf{H}_n^{-1/2} - \mathbf{I} = o_P(1)$, by Theorem 2(i) and Proposition 4. The result follows since $\|\mathbf{M}_n^{-1/2}\mathbf{H}_n^{1/2}\| \leq c_n^{1/2}$ and $\|\mathbf{H}_n^{1/2}(\hat{\beta}_n - \beta_0)\| \leq (\tilde{\tau}_n)^{1/2}r$.
□

Let $\gamma_n^{(D)} := \max_{1 \leq i \leq n} \lambda_{\max}(\mathbf{H}_n^{-1/2} \mathbf{X}_i^T \mathbf{A}_i^{1/2} \mathbf{R}_n^{-1} \mathbf{A}_i^{1/2} \mathbf{X}_i \mathbf{H}_n^{-1/2})$. Note that $\gamma_n^{(D)} \leq C d_n \tilde{\gamma}_n$, where $d_n = \max_{i \leq n, j \leq m} \sigma_{ij}^2$. We consider the following conditions:

$(\widetilde{\mathrm{N}}_\delta)$ there exists a $\delta > 0$ such that:
    (i) $Y := \sup_{i \geq 1} E(\|\mathbf{y}_i^*\|^{2+\delta} | \mathcal{F}_{i-1}) < \infty$ a.s.;
    (ii) $(c_n \tilde{\tau}_n)^{1+2/\delta} \gamma_n^{(D)} \to 0$,
$(\mathrm{C2})'$ $\max_{i \leq n} \lambda_{\max}(\mathcal{V}_i) \xrightarrow{P} 0$.

REMARK 6. Note that condition $(\widetilde{\mathrm{N}}_\delta)(\mathrm{i})$, with $Y$ integrable, implies condition (C1), whereas condition $(\mathrm{C2})'$ is a stronger form of (C2). Part (ii) of condition $(\widetilde{\mathrm{N}}_\delta)$ was introduced by Xie and Yang (2003).

The following result gives the asymptotic distribution of $\tilde{\mathbf{g}}_n$.

LEMMA 1. *Suppose that the conditions of Theorem* 1 *hold. Under conditions* $(\widetilde{\mathrm{N}}_\delta)$, $(\mathrm{C2})'$ *and* $(\mathrm{H}')$

$$\mathbf{M}_n^{-1/2} \tilde{\mathbf{g}}_n \xrightarrow{d} N(0, \mathbf{I}).$$

PROOF. We note that

$$\mathbf{M}_n^{-1/2} \tilde{\mathbf{g}}_n = \mathbf{M}_n^{-1/2} \mathbf{g}_n + \mathbf{M}_n^{-1/2} (\tilde{\mathbf{g}}_n - \mathbf{g}_n)$$

and $\|\mathbf{M}_n^{-1/2}(\tilde{\mathbf{g}}_n - \mathbf{g}_n)\| \leq (c_n \tilde{\tau}_n)^{1/2} \|(\tilde{\tau}_n)^{-1/2} \mathbf{H}_n^{-1/2} (\tilde{\mathbf{g}}_n - \mathbf{g}_n)\| \xrightarrow{P} 0$, by Proposition 3. Therefore it is enough to prove that $\mathbf{M}_n^{-1/2} \mathbf{g}_n \xrightarrow{d} N(0, \mathbf{I})$. By the Cramér–Wold theorem, this is equivalent to showing that: $\forall \lambda, \|\lambda\| = 1$

$$\lambda^T \mathbf{M}_n^{-1/2} \mathbf{g}_n = \sum_{i=1}^n Z_{n,i} \xrightarrow{d} N(0, 1), \tag{9}$$

where $Z_{n,i} = \lambda^T \mathbf{M}_n^{-1/2} \mathbf{X}_i^T \mathbf{A}_i^{1/2} \mathbf{R}_n^{-1} \mathbf{A}_i^{-1/2} \varepsilon_i$. Note that $E(Z_{n,i} | \mathcal{F}_{i-1}) = 0$ for all $i \leq n$, that is, $\{Z_{n,i}; i \leq n, n \geq 1\}$ is a martingale difference array.

Relationship (9) follows by the martingale central limit theorem with the Lindeberg condition [see Corollary 3.1 of Hall and Heyde (1980)] if

$$\sum_{i=1}^n E[Z_{n,i}^2 \mathbf{I}(|Z_{n,i}| > \varepsilon) | \mathcal{F}_{i-1}] \to 0 \quad \text{a.s.} \tag{10}$$

and

$$\sum_{i=1}^n E(Z_{n,i}^2 | \mathcal{F}_{i-1}) \xrightarrow{P} 1. \tag{11}$$

Relationship (10) follows from condition $(\widetilde{N}_\delta)$ exactly as in Lemma 2 of Xie and Yang (2003) with $\psi(t) = t^{\delta/2}$. Relationship (11) follows from conditions (C2)' and (H'):

$$\sum_{i=1}^n E(Z_{n,i}^2 | \mathcal{F}_{i-1}) - 1$$

$$= \sum_{i=1}^n [E(Z_{n,i}^2 | \mathcal{F}_{i-1}) - E(Z_{n,i}^2)]$$

$$= \sum_{i=1}^n \lambda^T \mathbf{M}_n^{-1/2} \mathbf{X}_i^T \mathbf{A}_i^{1/2} \mathbf{R}_n^{-1} \mathcal{V}_i \mathbf{R}_n^{-1} \mathbf{A}_i^{1/2} \mathbf{X}_i \mathbf{M}_n^{-1/2} \lambda$$

$$\leq \max_{1 \leq i \leq n} \lambda_{\max}(\mathcal{V}_i) \cdot \max_{1 \leq i \leq n} \lambda_{\max}(\bar{\mathbf{R}}_i^{-1}) \cdot \lambda^T \mathbf{M}_n^{-1/2} \mathbf{M}_n \mathbf{M}_n^{-1/2} \lambda$$

$$\leq C^{-1} \max_{i \leq n} \lambda_{\max}(\mathcal{V}_i) \xrightarrow{P} 0. \qquad \square$$

Putting together the results in Theorem 4 and Lemma 1, we obtain the asymptotic normality of the estimator $\hat{\beta}_n$.

THEOREM 5. *Under the conditions of Theorem 3 and conditions* $(\widetilde{N}_\delta)$, (C2)' *and* (H'),

$$\mathbf{M}_n^{-1/2} \mathbf{H}_n (\hat{\beta}_n - \beta_0) \xrightarrow{d} N(0, \mathbf{I}).$$

REMARK 7. In applications we would need a version of Theorem 5 where $\mathbf{M}_n$ is replaced by a consistent estimator. We suggest the estimator proposed by Liang and Zeger (1986) [see also Remark 8 of Xie and Yang (2003)]. The details of the proof are omitted.

## APPENDIX

**A.1.** The following lemma is a consequence of Kaufmann's (1987) martingale strong law of large numbers and can be viewed as a stronger version of Theorem 2.19 of Hall and Heyde (1980).

LEMMA A.1. *Let $(x_i)_{i \geq 1}$ be a sequence of random variables and let $(\mathcal{F}_i)_{i \geq 1}$ be a sequence of increasing $\sigma$-fields such that $x_i$ is $\mathcal{F}_i$-measurable for every $i \geq 1$. Suppose that $\sup_i E|x_i|^\alpha < \infty$ for some $\alpha \in (1, 2]$. Then*

$$\frac{1}{n} \sum_{i=1}^n (x_i - E(x_i | \mathcal{F}_{i-1})) \to 0 \qquad a.s. \text{ and in } L^\alpha.$$

PROOF. Note that $y_i = x_i - E(x_i|\mathcal{F}_{i-1})$, $n \geq 1$, is a martingale difference sequence. By the conditional Jensen inequality

$$|y_i|^\alpha \leq 2^{\alpha-1}\{|x_i|^\alpha + |E(x_i|\mathcal{F}_{i-1})|^\alpha\} \leq 2^{\alpha-1}\{|x_i|^\alpha + E(|x_i|^\alpha|\mathcal{F}_{i-1})\}$$

and $\sup_{i\geq 1} E|y_i|^\alpha \leq 2^\alpha \sup_{i\geq 1} E|x_i|^\alpha < \infty$. Hence

$$\sum_{i\geq 1} \frac{E|y_i|^\alpha}{i^\alpha} \leq \sup_{i\geq 1} E|y_i|^\alpha \cdot \sum_{i\geq 1} \frac{1}{i^\alpha} < \infty.$$

The lemma follows by Theorem 2 of Kaufmann (1987) with $p=1, B_i = i^{-1}$. □

PROOF OF PROPOSITION 1. We denote by $\hat{r}_{jk}^{(n)}, \bar{\bar{r}}_{jk}^{(n)}, v_{jk}^{(n)}$ ($j,k=1,\ldots,m$) the elements of the matrices $\widehat{\mathcal{R}}_n, \bar{\bar{\mathbf{R}}}_n, \mathcal{V}_n$, respectively. We write

$$(12) \quad \hat{r}_{jk}^{(n)} - \bar{\bar{r}}_{jk}^{(n)} = \frac{1}{n}\sum_{i=1}^n (y_{ij}^* y_{ik}^* - E(y_{ij}^* y_{ik}^*|\mathcal{F}_{i-1})) + \frac{1}{n}\sum_{i=1}^n v_{jk}^{(i)}.$$

The first term converges to zero almost surely and in $L^{1+\delta/2}$ by applying Lemma A.1 with $x_i = y_{ij}^* y_{ik}^*$, and using condition (C1). The second term converges to zero in probability by condition (C2). This convergence is also in $L^{1+\delta/2}$ because the sequence $\{n^{-1}\sum_{i=1}^n v_{jk}^{(i)}\}_n$ has uniformly bounded moments of order $1+\delta/2$ and hence is uniformly integrable. □

PROOF OF PROPOSITION 2. We denote by $\tilde{r}_{jk}^{(n)}$ ($j,k=1,\ldots,m$) the elements of the matrix $\widetilde{\mathcal{R}}_n$. Let $\tilde{\delta}_{i,jk} := [\sigma_{ij}\sigma_{ik}]/[\sigma_{ij}(\tilde{\beta}_n)\sigma_{ik}(\tilde{\beta}_n)] - 1$, $\widetilde{\Delta}\mu_{ij} := \mu_{ij}(\tilde{\beta}_n) - \mu_{ij}(\beta_0)$ and

$$\widetilde{\Delta}(\varepsilon_{ij}\varepsilon_{ik}) := \varepsilon_{ij}(\tilde{\beta}_n)\varepsilon_{ik}(\tilde{\beta}_n) - \varepsilon_{ij}\varepsilon_{ik} = (\widetilde{\Delta}\mu_{ij})(\widetilde{\Delta}\mu_{ik}) - (\widetilde{\Delta}\mu_{ij})\varepsilon_{ik} - (\widetilde{\Delta}\mu_{ik})\varepsilon_{ij}.$$

With this notation, we have

$$\tilde{r}_{jk}^{(n)} - \hat{r}_{jk}^{(n)} = \frac{1}{n}\sum_{i=1}^n \frac{\varepsilon_{ij}(\tilde{\beta}_n)\varepsilon_{ik}(\tilde{\beta}_n)}{\sigma_{ij}(\tilde{\beta}_n)\sigma_{ik}(\tilde{\beta}_n)} - \frac{1}{n}\sum_{i=1}^n \frac{\varepsilon_{ij}\varepsilon_{ik}}{\sigma_{ij}\sigma_{ik}}$$

$$= \frac{1}{n}\sum_{i=1}^n \frac{\widetilde{\Delta}(\varepsilon_{ij}\varepsilon_{ik})}{\sigma_{ij}\sigma_{ik}} + \frac{1}{n}\sum_{i=1}^n \frac{\widetilde{\Delta}(\varepsilon_{ij}\varepsilon_{ik})}{\sigma_{ij}\sigma_{ik}}\tilde{\delta}_{i,jk} + \frac{1}{n}\sum_{i=1}^n \frac{\varepsilon_{ij}\varepsilon_{ik}}{\sigma_{ij}\sigma_{ik}}\tilde{\delta}_{i,jk}.$$

From here, we conclude that

$$|\tilde{r}_{jk}^{(n)} - \hat{r}_{jk}^{(n)}| \leq U_{n,jk} + \max_{i\leq n}|\tilde{\delta}_{i,jk}| \cdot \left\{U_{n,jk} + \frac{1}{n}\sum_{i=1}^n |y_{ij}^* y_{ik}^*|\right\},$$

where

$$U_{n,jk} := \frac{1}{n}\sum_{i=1}^{n} \frac{|\widetilde{\Delta}\mu_{ij}| \cdot |\widetilde{\Delta}\mu_{ik}|}{\sigma_{ij}\sigma_{ik}} + \frac{1}{n}\sum_{i=1}^{n} \frac{|\widetilde{\Delta}\mu_{ij}|}{\sigma_{ij}} \cdot |y_{ik}^*| + \frac{1}{n}\sum_{i=1}^{n} \frac{|\widetilde{\Delta}\mu_{ik}|}{\sigma_{ik}} \cdot |y_{ij}^*|$$

$$= U_{n,jk}^{[1]} + U_{n,jk}^{[2]} + U_{n,jk}^{[3]}.$$

Recall that our estimator $\tilde{\beta}_n$ was obtained in the proof of Theorem 2 of Xie and Yang (2003) as a solution of the GEE in the case when all the "working" correlation matrices are $\mathbf{R}_i^{\mathrm{indep}} = \mathbf{I}$. One of the consequences of the result of Xie and Yang is that for every fixed $\varepsilon > 0$, there exist $r = r_\varepsilon$ and $N = N_\varepsilon$ such that, if we denote $\Omega_{n,\varepsilon} = \{\tilde{\beta}_n \text{ lies in } B_n^{\mathrm{indep}}(r)\}$, then

$$P(\Omega_{n,\varepsilon}) \geq 1 - \varepsilon \qquad \text{for all } n \geq N.$$

We define $\tilde{\beta}_n$ to be equal to $\beta_0$ on the event $\Omega_{n,\varepsilon}^c$. Therefore,

$$\text{on } \Omega_{n,\varepsilon}^c : \max_{i \leq n} |\tilde{\delta}_{i,jk}| = 0 \quad \text{and} \quad \widetilde{\Delta}\mu_{ij} = 0.$$

Using Taylor's formula and condition (AH)$^{\mathrm{indep}}$, we can conclude that on the event $\Omega_{n,\varepsilon}$, there exists a constant $C = C_\varepsilon$ such that

$$|\tilde{\delta}_{i,jk}| = \left|\frac{\dot{\mu}(\mathbf{x}_{ij}^T \beta_0)}{\dot{\mu}(\mathbf{x}_{ij}^T \tilde{\beta}_n)} - 1\right| \leq C \cdot (\gamma_n^0)^{\mathrm{indep}} \cdot (m^{1/2} r) \qquad \text{for all } i \leq n,$$

$$\frac{1}{n}\sum_{i=1}^{n} \frac{(\widetilde{\Delta}\mu_{ij})^2}{\sigma_{ij}^2} = n^{-1}(\tilde{\beta}_n - \beta_0)^T \left\{\sum_{i=1}^{n} \left(\frac{\sigma_{ij}^2(\bar{\beta}_n)}{\sigma_{ij}^2}\right)^2 \sigma_{ij}^2 \mathbf{x}_{ij}\mathbf{x}_{ij}^T\right\}(\tilde{\beta}_n - \beta_0)$$

$$\leq n^{-1}(\tilde{\beta}_n - \beta_0)^T \left\{\sum_{i=1}^{n} \mathbf{X}_i^T \mathbf{A}_i^{1/2} [\mathbf{A}_i(\bar{\beta}_n)\mathbf{A}_i^{-1}]^2 \mathbf{A}_i^{1/2} \mathbf{X}_i\right\}(\tilde{\beta}_n - \beta_0)$$

$$\leq n^{-1} \max_{i \leq n} \lambda_{\max}^2[\mathbf{A}_i(\bar{\beta}_n)\mathbf{A}_i^{-1}] \cdot \|(\mathbf{H}_n^{\mathrm{indep}})^{1/2}(\tilde{\beta}_n - \beta_0)\|^2$$

$$\leq C n^{-1}(m^2 r).$$

Note also that $E[n^{-1}\sum_{i=1}^{n}(y_{ij}^*)^2] = E[\hat{r}_{jj}^{(n)}] = O(1)$ since $\hat{r}_{jj}^{(n)} - \bar{\bar{r}}_{jj}^{(n)} \xrightarrow{L^1} 0$ and $\bar{\bar{r}}_{jj}^{(n)} = 1$. Applying the Cauchy–Schwarz inequality to each of the three sums that form $U_{n,jk}$, we can conclude that

$$E[U_{n,jk}^{[1]}] \to 0 \quad \text{and} \quad E[(U_{n,jk}^{[l]})^2] \to 0, \qquad l = 2, 3.$$

On the other hand,

$$E\left[\max_{i \leq n} |\tilde{\delta}_{i,jk}| \cdot U_{n,jk}\right] = \int_{\Omega_{n,\varepsilon}} \max_{i \leq n} |\tilde{\delta}_{i,jk}| \cdot U_{n,jk}\, dP$$

$$\leq C(\gamma_n^{(0)})^{\text{indep}} \int_\Omega U_{n,jk} \to 0,$$

$$E\left[\max_{i\leq n}|\tilde{\delta}_{i,jk}| \cdot \frac{1}{n}\sum_{i=1}^n |y_{ij}^* y_{ik}^*|\right] \leq C(\gamma_n^{(0)})^{\text{indep}} [E(\hat{r}_{jj}^{(n)})]^{1/2} [E(\hat{r}_{kk}^{(n)})]^{1/2} \to 0.$$

□

**A.2.**

PROOF OF PROPOSITION 3. Let $h_{ijk}^{(0)} = [\sigma_{ij}^2/\sigma_{ik}^2]^{1/2}$, $\widetilde{\mathcal{R}}_n^{-1} := \widetilde{\mathcal{Q}}_n = (\tilde{q}_{jk}^{(n)})_{j,k=1,\ldots,m}$ and $\mathbf{R}_n^{-1} := \mathbf{Q}_n = (q_{jk}^{(n)})_{j,k=1,\ldots,m}$. With this notation, we write

$$(\tilde{\tau}_n)^{-1/2} \mathbf{H}_n^{-1/2}(\tilde{\mathbf{g}}_n - \mathbf{g}_n)$$
$$= \sum_{j,k=1}^m (\tilde{q}_{jk}^{(n)} - q_{jk}^{(n)}) \cdot \left\{(\tilde{\tau}_n)^{-1/2} \mathbf{H}_n^{-1/2} \sum_{i=1}^n h_{ijk}^{(0)} \mathbf{x}_{ij} \varepsilon_{ik}\right\}.$$

By Theorem 1, $\tilde{q}_{jk}^{(n)} - q_{jk}^{(n)} \xrightarrow{P} 0$ for every $j,k$. The result will follow once we prove that $\{(\tilde{\tau}_n)^{-1/2}\mathbf{H}_n^{-1/2}\sum_{i=1}^n h_{ijk}^{(0)}\mathbf{x}_{ij}\varepsilon_{ik}\}_n$ is bounded in $L^2$ for every $j,k$. Since $(\varepsilon_{ik})_{i\geq 1}$ is a martingale difference sequence, we have

$$E\left(\left\|(\tilde{\tau}_n)^{-1/2}\mathbf{H}_n^{-1/2}\sum_{i=1}^n h_{ijk}^{(0)}\mathbf{x}_{ij}\varepsilon_{ik}\right\|^2\right)$$
$$= (\tilde{\tau}_n)^{-1} \operatorname{tr}\left\{\mathbf{H}_n^{-1/2}\left(\sum_{i=1}^n (h_{ijk}^{(0)})^2 \sigma_{ik}^2 \mathbf{x}_{ij}\mathbf{x}_{ij}^T\right)\mathbf{H}_n^{-1/2}\right\}$$
$$= (\tilde{\tau}_n)^{-1} \operatorname{tr}\left\{\mathbf{H}_n^{-1/2}\left(\sum_{i=1}^n \sigma_{ij}^2 \mathbf{x}_{ij}\mathbf{x}_{ij}^T\right)\mathbf{H}_n^{-1/2}\right\}$$
$$\leq (\tilde{\tau}_n)^{-1}(4m\tilde{\tau}_n)\operatorname{tr}(\mathbf{I}) = 4mp$$

because $\sum_{i=1}^n \sigma_{ij}^2 \mathbf{x}_{ij}\mathbf{x}_{ij}^T \leq \sum_{i=1}^n \mathbf{X}_i^T \mathbf{A}_i \mathbf{X}_i \leq \lambda_{\max}(\mathbf{R}_n)\mathbf{H}_n \leq 4m\tilde{\tau}_n \mathbf{H}_n$. □

PROOF OF PROPOSITION 4. We write

$$\mathcal{D}_n(\beta) = \mathbf{H}_n(\beta) + \mathbf{B}_n(\beta) + \mathcal{E}_n(\beta), \qquad \widetilde{\mathcal{D}}_n(\beta) = \widetilde{\mathcal{H}}_n(\beta) + \widetilde{\mathcal{B}}_n(\beta) + \widetilde{\mathcal{E}}_n(\beta),$$

where $\widetilde{\mathcal{H}}_n(\beta), \widetilde{\mathcal{B}}_n(\beta), \widetilde{\mathcal{E}}_n(\beta)$ have the same expressions as $\mathbf{H}_n(\beta), \mathbf{B}_n(\beta), \mathcal{E}_n(\beta)$, with $\mathbf{R}_n$ replaced by $\widetilde{\mathcal{R}}_n$. Our result will follow by the following three lemmas.

□

LEMMA A.2. *Suppose that condition* $(\widetilde{\mathrm{AH}})$ *holds. If* $(\pi_n \tilde{\gamma}_n)_n$ *is bounded, then for any* $r > 0$ *and for any* $p \times 1$ *vector* $\lambda$ *with* $\|\lambda\| = 1$,

$$\sup_{\beta \in \widetilde{B}_n(r)} |\lambda^T \mathbf{H}_n^{-1/2}[\widetilde{\mathcal{H}}_n(\beta) - \mathbf{H}_n(\beta)]\mathbf{H}_n^{-1/2}\lambda| \xrightarrow{P} 0.$$

LEMMA A.3. *Suppose that condition* $(\widetilde{\mathrm{AH}})$ *holds. If* $(\pi_n^2 \tilde{\gamma}_n)_n$ *is bounded, then for any* $r > 0$ *and for any* $p \times 1$ *vector* $\lambda$ *with* $\|\lambda\| = 1$,

$$\sup_{\beta \in \widetilde{B}_n(r)} |\lambda^T \mathbf{H}_n^{-1/2}[\widetilde{\mathcal{B}}_n(\beta) - \mathbf{B}_n(\beta)]\mathbf{H}_n^{-1/2}\lambda| \xrightarrow{P} 0.$$

LEMMA A.4. *Suppose that condition* $(\widetilde{\mathrm{AH}})$ *holds. If* $(n^{1/2} \tilde{\gamma}_n)_n$ *is bounded, then for any* $r > 0$ *and for any* $p \times 1$ *vector* $\lambda$ *with* $\|\lambda\| = 1$,

$$\sup_{\beta \in \widetilde{B}_n(r)} |\lambda^T \mathbf{H}_n^{-1/2}[\widetilde{\mathcal{E}}_n(\beta) - \mathcal{E}_n(\beta)]\mathbf{H}_n^{-1/2}\lambda| \xrightarrow{P} 0.$$

PROOF OF LEMMA A.2. Using Theorem 1 and the fact that $|r_{jk}^{(n)}| \leq 2$ for $n$ large, we have $\mathcal{A}_n = \mathbf{R}_n^{1/2}\widetilde{\mathcal{R}}_n^{-1}\mathbf{R}_n^{1/2} - \mathbf{I} = \mathbf{R}_n^{1/2}(\widetilde{\mathcal{R}}_n^{-1} - \mathbf{R}_n^{-1})\mathbf{R}_n^{1/2} \xrightarrow{P} 0$ (elementwise). For every $\beta$,

$$|\lambda^T \mathbf{H}_n^{-1/2}[\widetilde{\mathcal{H}}_n(\beta) - \mathbf{H}_n(\beta)]\mathbf{H}_n^{-1/2}\lambda|$$
$$= \left|\sum_{i=1}^n \lambda^T \mathbf{H}_n^{-1/2}\mathbf{X}_i^T \mathbf{A}_i(\beta)^{1/2}\mathbf{R}_n^{-1/2}\mathcal{A}_n \mathbf{R}_n^{-1/2}\mathbf{A}_i(\beta)^{1/2}\mathbf{X}_i\mathbf{H}_n^{-1/2}\lambda\right|$$
$$\leq \max\{|\lambda_{\max}(\mathcal{A}_n)|, |\lambda_{\min}(\mathcal{A}_n)|\} \cdot \{\lambda^T \mathbf{H}_n^{-1/2}\mathbf{H}_n(\beta)\mathbf{H}_n^{-1/2}\lambda\}.$$

The result follows, since one can show that for every $\beta \in \widetilde{B}_n(r)$

$$|\lambda^T \mathbf{H}_n^{-1/2}\mathbf{H}_n(\beta)\mathbf{H}_n^{-1/2}\lambda - 1|$$
(13)
$$\leq \lambda^T \mathbf{H}_n^{-1/2}\mathbf{H}_n^{[1]}(\beta)\mathbf{H}_n^{-1/2}\lambda + 2|\lambda^T \mathbf{H}_n^{-1/2}\mathbf{H}_n^{[2]}(\beta)\mathbf{H}_n^{-1/2}\lambda|$$
$$\leq C\pi_n\tilde{\gamma}_n + 2C(\pi_n\tilde{\gamma}_n)^{1/2} \leq C,$$

where

$$\mathbf{H}_n^{[1]}(\beta) = \sum_{i=1}^n \mathbf{X}_i^T(\mathbf{A}_i^{1/2}(\beta) - \mathbf{A}_i^{1/2})\mathbf{R}_n^{-1}(\mathbf{A}_i^{1/2}(\beta) - \mathbf{A}_i^{1/2})\mathbf{X}_i,$$

$$\mathbf{H}_n^{[2]}(\beta) = \sum_{i=1}^n \mathbf{X}_i^T(\mathbf{A}_i^{1/2}(\beta) - \mathbf{A}_i^{1/2})\mathbf{R}_n^{-1}\mathbf{A}_i^{1/2}\mathbf{X}_i.$$

We used the fact that $\sup_{\beta \in \widetilde{B}_n(r)} \max_{i \leq n} \lambda_{\max}\{(\mathbf{A}_i^{1/2}(\beta)\mathbf{A}_i^{-1/2} - \mathbf{I})^2\} \leq C\widetilde{\gamma}_n$, which follows by condition (A$\widetilde{\text{H}}$) as in Lemma B.1(ii) of Xie and Yang (2003). $\square$

PROOF OF LEMMA A.3. Let $\mathbf{w}_{i,n}(\beta)^T = \lambda^T \mathbf{H}_n^{-1/2} \mathbf{X}_i^T \mathbf{G}_i^{[1]}(\beta) \times \operatorname{diag}\{\mathbf{X}_i \mathbf{H}_n^{-1/2}\lambda\} \mathbf{R}_n^{-1/2}$ and $\mathbf{z}_{i,n}(\beta) = \mathbf{R}_n^{-1/2} \mathbf{A}_i(\beta)^{-1/2}(\mu_i - \mu_i(\beta))$. We have

$$|\lambda^T \mathbf{H}_n^{-1/2}[\widetilde{\mathcal{B}}_n^{[1]}(\beta) - \mathbf{B}_n^{[1]}(\beta)]\mathbf{H}_n^{-1/2}\lambda|$$
$$= \left|\sum_{i=1}^n \mathbf{w}_{i,n}(\beta)^T \mathcal{A}_n \mathbf{z}_{i,n}(\beta)\right|$$
$$\leq \|\mathcal{A}_n\| \left\{\sum_{i=1}^n \|\mathbf{w}_{i,n}(\beta)\|^2\right\}^{1/2} \left\{\sum_{i=1}^n \|\mathbf{z}_{i,n}(\beta)\|^2\right\}^{1/2}$$

by using the Cauchy–Schwarz inequality. Methods similar to those developed in the proof of Lemma A.2(ii) of Xie and Yang (2003) show that for any $\beta \in \widetilde{B}_n(r)$, $\sum_{i=1}^n \|\mathbf{w}_{i,n}(\beta)\|^2 \leq C\pi_n \gamma_n^{(0)}$ and $\sum_{i=1}^n \|\mathbf{z}_{i,n}(\beta)\|^2 \leq C\pi_n \widetilde{\tau}_n \lambda_{\max}(\mathbf{H}_n^{-1/2}\mathbf{H}_n(\bar{\beta}) \times \mathbf{H}_n^{-1/2}) \leq C\pi_n \widetilde{\tau}_n$ [using (13) for the last inequality]. Hence

$$\sup_{\beta \in \widetilde{B}_n(r)} |\lambda^T \mathbf{H}_n^{-1/2}[\widetilde{\mathcal{B}}_n^{[1]}(\beta) - \mathbf{B}_n^{[1]}(\beta)]\mathbf{H}_n^{-1/2}\lambda| \leq C\|\mathcal{A}_n\|\pi_n(\widetilde{\gamma}_n)^{1/2} \xrightarrow{P} 0.$$

Let $\mathbf{v}_{i,n}(\beta)^T = \lambda^T \mathbf{H}_n^{-1/2} \mathbf{X}_i^T \mathbf{A}_i(\beta)^{1/2} \mathbf{R}_n^{-1/2} \mathcal{A}_n \mathbf{R}_n^{-1/2} \operatorname{diag}\{\mathbf{X}_i \mathbf{H}_n^{-1/2}\lambda\} \mathbf{G}_i^{[2]}(\beta) \times \mathbf{A}_i(\beta)^{1/2} \mathbf{R}_n^{1/2}$. We have

$$|\lambda^T \mathbf{H}_n^{-1/2}[\widetilde{\mathcal{B}}_n^{[2]}(\beta) - \mathbf{B}_n^{[2]}(\beta)]\mathbf{H}_n^{-1/2}\lambda|$$
$$= \left|\sum_{i=1}^n \mathbf{v}_{i,n}(\beta)^T \mathbf{z}_{i,n}(\beta)\right| \leq \left\{\sum_{i=1}^n \|\mathbf{v}_{i,n}(\beta)\|^2\right\}^{1/2} \left\{\sum_{i=1}^n \|\mathbf{z}_{i,n}(\beta)\|^2\right\}^{1/2}.$$

One can prove that for any $\beta \in \widetilde{B}_n(r)$, $\sum_{i=1}^n \|\mathbf{v}_{i,n}(\beta)\|^2 \leq C\pi_n \gamma_n^{(0)} \lambda_{\max}(\mathcal{A}_n^2) \times \{\lambda^T \mathbf{H}_n^{-1/2}\mathbf{H}_n(\beta)\mathbf{H}_n^{-1/2}\lambda\} \leq C\pi_n \gamma_n^{(0)} \|\mathcal{A}_n\|^2$ [using (13) for the last inequality]. Hence

$$\sup_{\beta \in \widetilde{B}_n(r)} |\lambda^T \mathbf{H}_n^{-1/2}[\widetilde{\mathcal{B}}_n^{[2]}(\beta) - \mathbf{B}_n^{[2]}(\beta)]\mathbf{H}_n^{-1/2}\lambda| \leq C\|\mathcal{A}_n\|\pi_n(\widetilde{\gamma}_n)^{1/2} \xrightarrow{P} 0. \quad \square$$

PROOF OF LEMMA A.4. We write $\widetilde{\mathcal{E}}_n(\beta) - \mathcal{E}_n(\beta) = [\widetilde{\mathcal{E}}_n^{[1]}(\beta) - \mathcal{E}_n^{[1]}(\beta)] + [\widetilde{\mathcal{E}}_n^{[2]}(\beta) - \mathcal{E}_n^{[2]}(\beta)]$ and we use a decomposition which is similar to that given in the proof of Lemma A.3(ii) of Xie and Yang (2003). More precisely, we write

$$\lambda^T \mathbf{H}_n^{-1/2}[\widetilde{\mathcal{E}}_n^{[1]}(\beta) - \mathcal{E}_n^{[1]}(\beta)]\mathbf{H}_n^{-1/2}\lambda = T_n^{[1]} + T_n^{[3]}(\beta) + T_n^{[5]}(\beta),$$
$$\lambda^T \mathbf{H}_n^{-1/2}[\widetilde{\mathcal{E}}_n^{[2]}(\beta) - \mathcal{E}_n^{[2]}(\beta)]\mathbf{H}_n^{-1/2}\lambda = T_n^{[2]} + T_n^{[4]}(\beta) + T_n^{[6]}(\beta),$$

ESTIMATION AND LONGITUDINAL DATA 19where $T_n^{[l]}(\beta) = \sum_{j,k=1}^m (\tilde{q}_{jk}^{(n)} - q_{jk}^{(n)}) \cdot S_{n,jk}^{[l]}(\beta)$ for $l = 1, \ldots, 6$ and

$$S_{n,jk}^{[1]} = \lambda^T \mathbf{H}_n^{-1/2} \sum_{i=1}^n [\mathbf{A}_i^{-1/2} \mathbf{G}_i^{[1]}]_j [\operatorname{diag}\{\mathbf{X}_i \mathbf{H}_n^{-1/2} \lambda\}]_j h_{ijk}^{(0)} \mathbf{x}_{ij} \varepsilon_{ik},$$

$$S_{n,jk}^{[3]}(\beta) = \lambda^T \mathbf{H}_n^{-1/2} \sum_{i=1}^n [\mathbf{A}_i^{-1/2} \mathbf{G}_i^{[1]}(\beta)]_j [\operatorname{diag}\{\mathbf{X}_i \mathbf{H}_n^{-1/2} \lambda\}]_j [\mathbf{A}_i(\beta)^{-1/2} \mathbf{A}_i^{1/2} - I]_k$$
$$\times h_{ijk}^{(0)} \mathbf{x}_{ij} \varepsilon_{ik},$$

$$S_{n,jk}^{[5]}(\beta) = \lambda^T \mathbf{H}_n^{-1/2} \sum_{i=1}^n [\mathbf{A}_i^{-1/2} (\mathbf{G}_i^{[1]}(\beta) - \mathbf{G}_i^{[1]})]_j [\operatorname{diag}\{\mathbf{X}_i \mathbf{H}_n^{-1/2} \lambda\}]_j h_{ijk}^{(0)} \mathbf{x}_{ij} \varepsilon_{ik},$$

$$S_{n,jk}^{[2]} = \lambda^T \mathbf{H}_n^{-1/2} \sum_{i=1}^n [\operatorname{diag}\{\mathbf{X}_i \mathbf{H}_n^{-1/2} \lambda\}]_k [\mathbf{G}_i^{[2]} \mathbf{A}_i^{1/2}]_k h_{ijk}^{(0)} \mathbf{x}_{ij} \varepsilon_{ik},$$

$$S_{n,jk}^{[4]}(\beta) = \lambda^T \mathbf{H}_n^{-1/2} \sum_{i=1}^n [\mathbf{A}_i^{-1/2} \mathbf{A}_i(\beta)^{1/2} - I]_j [\operatorname{diag}\{\mathbf{X}_i \mathbf{H}_n^{-1/2} \lambda\}]_k [\mathbf{G}_i^{[2]}(\beta) \mathbf{A}_i^{1/2}]_k$$
$$\times h_{ijk}^{(0)} \mathbf{x}_{ij} \varepsilon_{ik},$$

$$S_{n,jk}^{[6]}(\beta) = \lambda^T \mathbf{H}_n^{-1/2} \sum_{i=1}^n [\operatorname{diag}\{\mathbf{X}_i \mathbf{H}_n^{-1/2} \lambda\}]_k [(\mathbf{G}_i^{[2]}(\beta) - \mathbf{G}_i^{[2]}) \mathbf{A}_i^{1/2}]_k h_{ijk}^{(0)} \mathbf{x}_{ij} \varepsilon_{ik}$$

(here we have denoted with $[\Delta]_j$ the $j$th element on the diagonal of a matrix $\Delta$).

Since $\tilde{q}_{jk}^{(n)} - q_{jk}^{(n)} \xrightarrow{P} 0$, it is enough to prove that $\{S_{n,jk}^{[1]}\}_n$, $\{S_{n,jk}^{[2]}\}_n$ and $\{\sup_{\beta \in \widetilde{B}_n(r)} |S_{n,jk}^{[l]}(\beta)|\}_n$, $l = 3, 4, 5, 6$, are bounded in $L^2$ for every $j, k = 1, \ldots, m$.

We have

$$E(|S_{n,jk}^{[1]}|^2)$$
$$\leq \operatorname{tr}\left\{\mathbf{H}_n^{-1/2} \left(\sum_{i=1}^n [\mathbf{A}_i^{-1/2} \mathbf{G}_i^{[1]}]_j^2 [\operatorname{diag}\{\mathbf{X}_i \mathbf{H}_n^{-1/2} \lambda\}]_j^2 (h_{ijk}^{(0)})^2 \sigma_{ik}^2 \mathbf{x}_{ij} \mathbf{x}_{ij}^T\right) \mathbf{H}_n^{-1/2}\right\}$$
$$\leq C \gamma_n^{(0)} \operatorname{tr}\left\{\mathbf{H}_n^{-1/2} \left(\sum_{i=1}^n \sigma_{ij}^2 \mathbf{x}_{ij} \mathbf{x}_{ij}^T\right) \mathbf{H}_n^{-1/2}\right\} \leq C \gamma_n^{(0)} (4mp\tilde{\tau}_n) = C \tilde{\gamma}_n \leq C.$$

By the Cauchy–Schwarz inequality, for every $\beta \in \widetilde{B}_n(r)$,

$$|S_{n,jk}^{[3]}(\beta)|^2 \leq \left\{\sum_{i=1}^n [\mathbf{A}_i^{-1/2} \mathbf{G}_i^{[1]}(\beta)]_j^2 [\operatorname{diag}\{\mathbf{X}_i \mathbf{H}_n^{-1/2} \lambda\}]_j^2 [\mathbf{A}_i(\beta)^{-1/2} \mathbf{A}_i^{1/2} - I]_k^2\right\}$$

$$\times \left\{ \sum_{i=1}^{n} (h_{ijk}^{(0)})^2 \varepsilon_{ik}^2 (\lambda^T \mathbf{H}_n^{-1/2} \mathbf{x}_{ij})^2 \right\}$$

$$\leq C n \gamma_n^{(0)} \tilde{\gamma}_n \cdot \left\{ \lambda^T \mathbf{H}_n^{-1/2} \left( \sum_{i=1}^{n} (h_{ijk}^{(0)})^2 \varepsilon_{ik}^2 \mathbf{x}_{ij} \mathbf{x}_{ij}^T \right) \mathbf{H}_n^{-1/2} \lambda \right\}.$$

Hence $E(\sup_{\beta \in \widetilde{B}_n(r)} |S_{n,jk}^{[3]}(\beta)|^2) \leq C n \gamma_n^{(0)} \tilde{\gamma}_n \cdot \{\lambda^T \mathbf{H}_n^{-1/2} (\sum_{i=1}^{n} \sigma_{ij}^2 \mathbf{x}_{ij} \mathbf{x}_{ij}^T) \mathbf{H}_n^{-1/2} \times \lambda\} \leq C n \gamma_n^0 \tilde{\gamma}_n \cdot (4m\tilde{\tau}_n) \leq C n (\tilde{\gamma}_n)^2 \leq C$.

Similarly, by the Cauchy–Schwarz inequality, for every $\beta \in \widetilde{B}_n(r)$,

$$|S_{n,jk}^{[5]}(\beta)|^2 \leq \left\{ \sum_{i=1}^{n} [\mathbf{A}_i^{-1/2} (\mathbf{G}_i^{[1]}(\beta) - \mathbf{G}_i^{[1]})]_j^2 [\text{diag}\{\mathbf{X}_i \mathbf{H}_n^{-1/2} \lambda\}]_j^2 \right\}$$

$$\times \left\{ \sum_{i=1}^{n} (h_{ijk}^{(0)})^2 \varepsilon_{ik}^2 (\lambda^T \mathbf{H}_n^{-1/2} \mathbf{x}_{ij})^2 \right\}$$

$$\leq C n \gamma_n^{(0)} \tilde{\gamma}_n \cdot \left\{ \lambda^T \mathbf{H}_n^{-1/2} \left( \sum_{i=1}^{n} (h_{ijk}^{(0)})^2 \varepsilon_{ik}^2 \mathbf{x}_{ij} \mathbf{x}_{ij}^T \right) \mathbf{H}_n^{-1/2} \lambda \right\}$$

and $E(\sup_{\beta \in \widetilde{B}_n(r)} |S_{n,jk}^{[5]}(\beta)|^2) \leq C n (\tilde{\gamma}_n)^2 \leq C$.

The terms $S_{n,jk}^{[l]}(\beta)$, $l = 2, 4, 6$, can be treated by similar methods. $\square$

## REFERENCES

- Chen, K., Hu, I. and Ying, Z. (1999). Strong consistency of maximum quasi-likelihood estimators in generalized linear models with fixed and adaptive designs. *Ann. Statist.* **27** 1155–1163. MR1740117
- Fahrmeir, L. and Kaufmann, H. (1985). Consistency and asymptotic normality of the maximum likelihood estimator in generalized linear models. *Ann. Statist.* **13** 342–368. MR773172
- Fahrmeir, L. and Tutz, G. (1994). *Multivariate Statistical Modelling Based on Generalized Linear Models.* Springer, New York. MR1284203
- Hall, P. and Heyde, C. C. (1980). *Martingale Limit Theory and Its Application.* Academic Press, New York. MR624435
- Kaufmann, H. (1987). On the strong law of large numbers for multivariate martingales. *Stochastic Process. Appl.* **26** 73–85. MR917247
- Lai, T. L., Robbins, H. and Wei, C. Z. (1979). Strong consistency of least squares estimates in multiple regression. II. *J. Multivariate Anal.* **9** 343–361. MR548786
- Liang, K.-Y. and Zeger, S. L. (1986). Longitudinal data analysis using generalized linear models. *Biometrika* **73** 13–22. MR836430
- McCullagh, P. and Nelder, J. A. (1989). *Generalized Linear Models*, 2nd ed. Chapman and Hall, London. MR727836
- Rao, J. N. K. (1998). Marginal models for repeated observations: Inference with survey data. In *Proc. Section on Survey Research Methods* 76–82. Amer. Statist. Assoc., Alexandria, VA.

SCHIOPU-KRATINA, I. (2003). Asymptotic results for generalized estimating equations with data from complex surveys. *Rev. Roumaine Math. Pures Appl.* **48** 327–342. MR2038208

SCHOTT, J. R. (1997). *Matrix Analysis for Statistics.* Wiley, New York. MR1421574

SHAO, J. (1992). Asymptotic theory in generalized linear models with nuisance scale parameters. *Probab. Theory Related Fields* **91** 25–41. MR1142760

SHAO, J. (1999). *Mathematical Statistics.* Springer, New York. MR1670883

XIE, M. and YANG, Y. (2003). Asymptotics for generalized estimating equations with large cluster sizes. *Ann. Statist.* **31** 310–347. MR1962509

YUAN, K.-H. and JENNRICH, R. I. (1998). Asymptotics of estimating equations under natural conditions. *J. Multivariate Anal.* **65** 245–260. MR1625893

DEPARTMENT OF MATHEMATICS
AND STATISTICS
UNIVERSITY OF OTTAWA
OTTAWA, ONTARIO
CANADA K1N 6N5
E-MAIL: rbala348@science.uottawa.ca

STATISTICS CANADA
HSMD 16RHC
OTTAWA, ONTARIO
CANADA K1A 0T6
E-MAIL: ioana.schiopu-kratina@statcan.ca